\documentclass[11pt]{amsart}

\usepackage[margin=1.5in]{geometry}

\usepackage{amsmath}
\usepackage{amsfonts}
\usepackage{amssymb}
\usepackage{amsthm}
\usepackage{enumerate}
\usepackage{graphicx}
\usepackage{hyperref}
\usepackage{enumitem}
\usepackage{comment}
\usepackage{color}

\usepackage{verbatim}

\input xy
\xyoption{all}

\theoremstyle{plain}
\newtheorem{theorem}{Theorem}[section]

\newtheorem{lemma}[theorem]{Lemma}
\newtheorem{proposition}[theorem]{Proposition}
\newtheorem*{proposition*}{Proposition}

\theoremstyle{definition}
\newtheorem{definition}[theorem]{Definition}

\newtheorem{rem}[theorem]{Remark}
\newtheorem{conv}[theorem]{Convention}

\DeclareMathOperator{\aut}{Aut}
\DeclareMathOperator{\fix}{Fix}
\newcommand{\cal}[1]{\mathcal{#1}}
\newcommand{\bb}[1]{\mathbb{#1}}
\newcommand{\stab}{\mathrm{Stab}}

\newcommand{\im}{\mathrm{im}}

\newcommand{\hh}{\mathfrak{h}}
\renewcommand{\hom}{\mathrm{Hom}}
\newcommand{\pgl}{\mathrm{PGL}}
\newcommand{\Lin}{\mathrm{Lin}}
\newcommand{\lin}{\mathrm{lin}}
\def\Z{{\mathbf{Z}}}

\def\C{{\mathbf{C}}}

\makeatletter
\def\ps@pprintTitle{%
  \let\@oddhead\@empty
  \let\@evenhead\@empty
  \let\@oddfoot\@empty
  \let\@evenfoot\@oddfoot
}
\makeatother

\usepackage{marvosym}

\title{Smooth quotients of complex tori by finite groups}
\author{Robert Auffarth}
\address{R. Auffarth \\Departamento de Matem\'aticas, Facultad de
Ciencias, Universidad de Chile\\ Las Palmeras 3425, \~Nu\~noa, Santiago, Chile}
\email{rfauffar@uchile.cl}
\author{Giancarlo Lucchini Arteche}
\address{G. Lucchini Arteche \\Departamento de Matem\'aticas, Facultad de
Ciencias, Universidad de Chile\\ Las Palmeras 3425, \~Nu\~noa, Santiago, Chile}
\email{luco@uchile.cl}

\thanks{The first author was partially supported by Conicyt via Fondecyt Grant 11180965. The second author was partially supported by Conicyt via Fondecyt Grant 11170016 and PAI Grant 79170034.}


\begin{document}

\begin{abstract}
Let $T$ be a complex torus and $G$ a finite group acting on $T$ without translations such that $T/G$ is smooth. Consider the subgroup $F\leq G$ generated by elements that have at least one fixed point. We prove that there exists a point $x\in T$ fixed by the whole group $F$ and that the quotient $T/G$ is a fibration of products of projective spaces over an \'etale quotient of a complex torus (the \'etale quotient being Galois with group $G/F$). In particular, when $G=F$, we may assume that $G$ fixes the origin. This is related to previous work by the authors, where the case of actions on abelian varieties fixing the origin was treated. Here, we generalize these results to complex tori and use them to reduce the problem of classifying smooth quotients of complex tori to the case of \'etale quotients.

An ingredient of the proof of our fixed-point theorem is a result proving that in every irreducible complex reflection group there is an element which is not contained in any proper reflection subgroup and that Coxeter elements have this property for well-generated groups. This result is proved by Stephen Griffeth in an appendix.

\end{abstract}

\maketitle

\section{Introduction}
Let $T$ be a complex torus of dimension $n$, and let $G$ be a finite subgroup of biholomorphisms of $T$. We are interested in studying when $T/G$ is smooth. This article is a twofold generalization of \cite{ALA}, where the same situation was considered, but for abelian varieties and in the case that $G$ fixes the origin. This article gives a more precise version of a result by Demailly, Hwang and Peternell \cite{Demailly}, which itself generalizes previous work on varieties covered by abelian varieties by Hwang and Mok \cite{HwangMok}. It can be seen as a bridge between our previous article and the theory of \textit{hyperelliptic manifolds}.

A hyperelliptic manifold is by definition a manifold $M$ that is not a complex torus, but such that there exists a complex torus $X$ and a finite subgroup $G\leq\mathrm{Aut}(X)$ that acts freely on $X$, such that $M\simeq X/G$. This generalizes the classic notion of hyperelliptic surface that was studied and characterized by Enriques and Severi \cite{ES} and Bagnera and de Franchis \cite{BF}. These manifolds have been studied in general by Lange \cite{Lange}, Catanese and Corvaja \cite{CC}, and have recently been characterized in dimension 3 by Catanese and Demleitner \cite{CatDem} who based their work on a characterization of discontinuous groups of affine transfor\-mations of $\C^3$ given by Uchida and Yoshihara \cite{UY}.

The purpose of our article is to reduce the classification of smooth quotients of complex tori  to the hyperelliptic case. We first note that we can assume that the group $G$ contains no translations, since if $\cal T_G$ is the (normal) subgroup of $G$ that consists of translations, then $G/\cal T_G$ acts on the complex torus $T/\cal T_G$ without translations. Secondly, let $F\leq G$ denote the subgroup of $G$ generated by elements that fix some point of $T$. We first note that this subgroup is clearly normal. Even though \emph{a priori} such a group could contain elements that do not fix any point, we prove the following (cf.~Section \ref{sec G=F}):

\begin{theorem}\label{fixedpointtheorem}
Let $G$ be a group acting faithfully on a complex torus $T$ with no translations and such that $T/G$ is smooth. Let $F$ be the (normal) subgroup generated by elements that fix at least a point. Then:
\begin{enumerate}
\item\label{FreeFaithful} $G/F$ acts freely on $T/F$;
\item\label{T/Fsmooth} $T/F$ is smooth;
\item\label{Ffixespoint} $F$ fixes a point $x\in T$.
\end{enumerate}
\end{theorem}

We prove further the following result (cf.~Section \ref{sec tori}), which provides a generalization of our results in \cite{ALA}:

\begin{proposition}\label{prop T0 and PG}
Let $F$ be a group that acts by biholomorphic homomorphisms on a complex torus $T$ such that $T/F$ is smooth and let $T_0$ be the connected component of $T^F$ containing 0. Then we have the following:
\begin{enumerate}
    \item\label{complementary} There exists a complementary $F$-stable subtorus $P_F$ of $T$ (i.e.~$T=T_0+P_F$ and $T_0\cap P_F$ is finite) which is an abelian variety. In particular, $T$ is an abelian variety if and only if $T_0$ is an abelian variety.
    \item\label{fibration} The quotient $P_F/F$ is isomorphic to a product of projective spaces and there exists a fibration $T/F\to T_0/(T_0\cap P_F)$ with fibers isomorphic to $P_F/F$.
\end{enumerate}
\end{proposition}

Thus, up to modifying the origin of $T$ in Theorem \ref{fixedpointtheorem}, we may apply these results to the pair $(T,F)$ and conclude that $T/F$ is a fibration of products of projective spaces over a complex torus $Z$ (which is isogenous to $T_0$). Since the group $G/F$ acts freely on this manifold, one may wonder what the quotient looks like. As it turns out, the action ``does not touch'' the fibers and acts directly on the complex torus $Z$. More precisely:

\begin{theorem}\label{thm etale quot}
Let $T,G,F$ be as above and let $\pi:T/F\to Z$ be the corresponding fibration. Then $G/F$ acts freely on $Z$, equivariantly with respect to $\pi$, and the quotient $T/G$ is isomorphic to a fibration of products of projective spaces over $Z/(G/F)$. In other words, the fibration $\pi$ is the pullback of a fibration over the quotient manifold $Z/(G/F)$ (which is either a complex torus or a hyperelliptic manifold).
\end{theorem}

This result, which we prove in Section \ref{sec reduction}, is a more precise version of \cite[Theorem 1.1 Part (3)]{Demailly}. Indeed, their result states that if $X$ is a holomorphic image of a complex torus, then $X$ is an \'etale quotient of the product of a complex torus with a product of projective spaces. The methods they use are very different from ours, and rely on a result by Hwang and Mok \cite[Main Theorem]{HwangMok} that states that if $X$ is projective it must be a projective bundle over another projective manifold (Demailly, Hwang and Peternell adapt this result so that it can be used in the case that $X$ is non-projective as well). This result is then used to obtain an \'etale trivialization of $X$ seen as a projective bundle. Their method has a ``from the ground up'' approach, in the sense that they start with $X$ and go up to the torus that is covering it in order to find the \'etale cover. In our scenario, on the other hand, we start off with a group acting on a complex torus and show explicitly in this case what the \'etale cover and Galois group are of the total quotient.

In particular, Theorem \ref{thm etale quot} implies that the full classification of smooth quotients of complex tori is now reduced to understanding fibrations of products of projective spaces over \'etale quotients of complex tori. We show in Section \ref{sec red hyperell case} how to obtain a classification of these fibrations when an \'etale quotient of a complex torus is fixed. More precisely, we show the exact way in which projective bundles can be obtained over a hyperelliptic manifold or complex torus, see Theorem \ref{thm classification}. Since the quotient of a complex torus by a finite group is an \'etale quotient of a trivial bundle over a complex torus, this quotient must be obtained by gluing products of projective spaces by automorphisms that descend from the complex torus. Using the classification results obtained in \cite{ALA}, Propositions \ref{prop class admissible irred} and \ref{prop class admissible general} give an explicit description of what possible ``gluing groups'' can appear.

It is worth noticing that, in the particular case where the base of the fibration is an abelian variety, a full classification of these gluing morphisms has been given by Mart\'inez-N\'u\~nez in \cite{Gary}. It is extremely likely that his results, like ours in \cite{ALA}, extend to the complex torus case without any modifications.

We finish the article with a section where we provide examples of non-abelian hyperelliptic manifolds of arbitrarily large dimension.

\subsection*{Acknowledgements}
We are indebted to Jean-Pierre Demailly, who pointed out to us the existence of \'etale quotients of complex tori and his result with Hwang and Peternell. We would also like to thank Sadek Al Harbat, Jean Michel and Don Taylor for helpful comments, and an anonymous referee for showing us literature on hyperelliptic manifolds. We thank as well a second anonymous referee, whose remarks helped us to improve the presentation of our results. And of course we thank Stephen Griffeth for taking the time to write the appendix to this article.

\section{Notations and Preliminaries}\label{sec notations}
Throughout the text, we will consider pairs $(T,G)$ consisting of a complex torus $T$ with an action of a finite group $G$. We will always assume that the action is faithful. As stated in the introduction, if we let $\cal T_G$ be the subgroup of $G$ consisting of translations of $T$, then it is an elementary exercise to show that $\cal T_G$ is normal in $G$, $G/\cal T_G$ acts on $T/\cal T_G$, and $T/G\cong (T/\cal T_G)/(G/\cal T_G)$. We thus make the following convention.

\begin{conv}
For the rest of the article, we consider pairs $(T,G)$ in which $G$ does not contain translations of the complex torus $T$.
\end{conv}

For any pair $(T,G)$, we will denote by $F$ the normal subgroup of $G$ generated by elements fixing at least a point. Note that $F$ is indeed normal since the conjugate of an element fixing a point also fixes a point. This need not imply that every element in $F$ fixes a point, but it will be the case \emph{a posteriori} in our setting.\\

We recall the definition of complex reflection groups and pseudoreflections, both in the classic setting and in our geometric context.

\begin{definition}
Let $G$ be an abstract group acting linearly on $\C^n$. We say that $g\in G$ is a \textit{pseudoreflection} if it fixes a hyperplane pointwise. We say that $G$ is a \textit{complex reflection group} if it is generated by pseudoreflections. Such a group is moreover called \textit{irreducible} if this is the case for the corresponding complex representation.
\end{definition}

Shephard and Todd gave a classification of complex reflection groups in \cite{ST}. As it turns out, these always split as direct products of irreducible ones. Moreover, irreducible complex reflection groups fall into three infinite families or 34 sporadic cases. Here we make an indirect use of this classification, since they are both used in the results that are found in \cite{ALA} as well as in the appendix to this article.

The main feature of complex reflection groups is that, given a linear action of a finite group $G$ on $\C^n$, the quotient variety $\C^n/G$ is smooth at the image of the origin if and only if $G$ is a complex reflection group. This can be generalized to arbitrary complex varieties by analyzing, for any point in the variety, the linear action of its stabilizer on its tangent space. This is why we need the following definition, which comes with an (innocuous) abuse of language.

\begin{definition}
Let $(T,G)$ be a complex torus $T$ with an action of a finite group $G$. Let $x\in T$ and $g\in G$. Then $g$ is a \textit{pseudoreflection} at $x$ if $\text{Fix}(g)$ is of codimension 1 and $x\in\fix(g)$.
\end{definition}

The Chevalley-Shephard-Todd Theorem states then that $T/G$ is smooth if and only if for every $x\in T$, $\stab_G(x)$ is a complex reflection group, i.e.~if it is generated by pseudoreflections (at $x$).\\

Recall that the group of biholomorphisms $\aut(T)$ has a natural structure of semi-direct product
\[\aut(T)=T\rtimes \Lin(T),\]
where $T$ corresponds to translations and $\Lin(T)$ denotes group automorphisms. The projection morphism $\pi:\aut(T)\to\Lin(T)$ is simply given by $\phi\mapsto \phi-\phi(0)$.

Given a pair $(T,G)$, we will denote by $G_\lin$ the image of $\pi:G\to\Lin(T)$. By our assumptions, we have that $G_\lin\cong G$ and $G_\lin$ fixes the origin. For an element $g\in G$, we will denote by $g_\lin$ its image $\pi(g)\in G_\lin$. We will abusively see $G_\lin$ acting at the same time on $T$ and on its universal cover $T_0(T)$ since it has a canonical lift.

\section{Proof of Theorem \ref{fixedpointtheorem}}\label{sec G=F}

For Item (\ref{FreeFaithful}) of Theorem \ref{fixedpointtheorem}, if the coset $gF$ has a fixed point on $T/F$ with a preimage $x\in T$, then in particular there exists $f\in F$ such that $g(x)=f(x)$. But then $f^{-1}g$ has a fixed point, and so $g\in F$. This proves the freeness of the action of $G/F$ on $T/F$.\\

For Item (\ref{T/Fsmooth}), by the previous paragraph we have that $T/F$ is an \'etale cover of $T/G$, and since the latter space is smooth we have that $T/F$ is smooth.\\

All we are left to prove then is the following result, which amounts to Theorem \ref{fixedpointtheorem} in the case $G=F$ and hence proves Item \eqref{Ffixespoint} of Theorem \ref{fixedpointtheorem}.

\begin{theorem}\label{fixedpointtheorem G=F}
Let $G$ be a group acting faithfully on a complex torus $T$ with no translations and such that $T/G$ is smooth. Assume that $G$ is generated by elements that fix at least a point. Then $G$ fixes a point $x\in T$.
\end{theorem}

This result can be reduced to the study of complex reflection groups as follows.

\begin{proposition}\label{prop G=F is crg}
Let $(T,G)$ be a pair as above. Then the action of $G_\lin$ on $T_0(T)$ realizes $G_\lin$ as a complex reflection group. Moreover, pseudoreflections in $G$ correspond to pseudoreflections in $G_\lin$ via the isomorphism $G\to G_\lin$.
\end{proposition}

\begin{proof}
Since $G$ is generated by elements that fix at least a point, we have
\[G=\langle\stab_G(x)\mid x\in T\rangle,\]
and since each of the subgroups in the equality is generated by pseudoreflections by the Chevalley-Shephard-Todd Theorem, we have that $G$ is generated by pseudo\-reflections (not all necessarily at the same point).

Now, let $g\in G$ be a pseudoreflection. Then its image $g_\lin\in G_\lin$ is a pseudoreflection as well. Indeed, since $g(x)=g_\lin(x)+g(0)$ for every $x\in T$, we see that
\[\fix(g)=(1-g_\lin)^{-1}(g(0))\qquad\text{and}\qquad \fix(g_\lin)=(1-g_\lin)^{-1}(0).\]
This means that $\fix(g)$ is the translate of $\fix(g_\lin)$ by any element $x\in\fix(g)$ (which is non-empty by assumption). In particular, since $\fix(g)$ has codimension 1, so does $\fix(g_\lin)$ and hence $g_\lin$ is a pseudoreflection fixing the origin. The same argument proves that $g$ is a pseudoreflection at some point $x$ if $g_\lin$ is a pseudoreflection.

Since $G$ is generated by pseudoreflections, we see that $G_\lin$ is generated by pseudo\-reflections fixing the origin. We obtain then that the analytic representation $\rho_a:G_\lin\to\mathrm{GL}(T_0(T))$ realizes $G_\lin\cong G$ as a complex reflection group.
\end{proof}

We can now use the theory of complex reflection groups to tackle Theorem \ref{fixedpointtheorem G=F}. As a general reference we recommend the classic articles by Shephard-Todd \cite{ST} and Cohen \cite{Cohen}. 

For starters, we know that $G_\lin\cong G_{1,\lin}\times\cdots\times G_{r,\lin}$ and $T_0(T)=W_0\oplus W_1\oplus\cdots\oplus W_r$ where 
\begin{itemize}
\item $W_i$ is an irreducible complex representation of $G_{i,\lin}$ that makes $G_{i,\lin}$ an irreducible finite complex reflection group for $i>0$;
\item $G_{j,\lin}$ acts trivially on $W_i$ for $i\neq j$.
\end{itemize}
In particular, $W_0=T_0(T)^{G_\lin}$.

We can then give a variant of \cite[Lem.~2.6]{ALA}. Write $G=G_1\times\cdots\times G_r$ so that $G_{i,\lin}$ is the image of $G_i$ in $\Lin(T)$.

\begin{lemma}\label{lemma desc Ti}
The subspace $W_i$ induces a $G_i$-stable (and $G_\lin$-stable) subtorus $T_i$ of $T$. In particular, the subtorus $P_G=\sum_{i=1}^r T_i$ is $G$-stable.
\end{lemma}

\begin{proof}
First of all, note that \cite[Lem.~2.6]{ALA} is completely valid in the context of complex tori since it does not use the projectivity assumption. This results tells us that the subtori $T_i$ are well defined, they are $G_\lin$-stable and $G_{i,\lin}$ acts trivially on $T_j$ for $i\neq j$. Let us prove that $T_i$ is $G_i$-stable as well. Every pseudoreflection $s\in G$ acts on $T$ as $s(x)=s_\lin(x)+z$ with $z\in T$. But from the proof of Proposition \ref{prop G=F is crg} above we deduce that $z\in\im(1-s_{0})$, otherwise $\fix(s)$ would be empty. Since $G_{i,\lin}$ acts trivially on $T_j$ for $j\neq i$, we see that if $s\in G_i$ then $\im(1-s_\lin)\subset T_i$, which proves the claim.

In particular, every pseudoreflection $s\in G$ acts on $T$ as $s(x)=s_\lin(x)+z$ with $z\in P_G$. Since $P_G$ is $G_\lin$-stable, we get that it is $G$-stable as well.
\end{proof}

Now we can even reduce the proof of Theorem \ref{fixedpointtheorem G=F} to a precise statement in the theory of complex reflection groups that is totally independent of our setting. This is done in the following result.

\begin{proposition}\label{prop fixed point}
Let $(T,G)$ be a pair as above. Assume moreover that there exists $g\in G$ such that
\begin{itemize}
\item $g_\lin$ is not contained in any proper complex reflection subgroup of $G_\lin$;
\item $g_\lin$ has no trivial eigenvalue.
\end{itemize}
Then $G$ fixes a point in $T$.
\end{proposition}

\begin{proof}
By the second assumption on $g_\lin$, we know that the endomorphism $(1-g_\lin)$ of $T_0(T)$ is invertible. This implies that the corresponding endomorphism of $T$ is surjective. Now, we know that $g$ acts on $T$ as $g(x)=g_\lin(x)+z$ for certain $z\in T$. Take $y\in T$ such that $(1-g_\lin)(y)=z$. Then clearly $y$ is fixed by $g$.

By the Chevalley-Shephard-Todd Theorem, we get then that the stabilizer $S\leq G$ of $y$ is generated by pseudoreflections and contains $g$. The first assumption implies then immediately that $S=G$, which concludes the proof.
\end{proof}

Note that Theorem \ref{fixedpointtheorem G=F} follows immediately from applying Proposition \ref{prop fixed point} to the pair $(P_G,G)$, which comes naturally from the pair $(T,G)$ thanks to Lemma \ref{lemma desc Ti}. Thus, all we have to do is to prove that, for every pair $(T,G)$ with $G_\lin$ acting on $T_0(T)=\C^n$ as a finite complex reflection group of rank $n$, we can find an element $g\in G$ such that $g_\lin$ is not contained in any proper complex reflection subgroup of $G_\lin$ and has no trivial eigenvalue. Now, since the existence of such a $g$ depends only on the action of $G_\lin$ on $T_0(T)$, we immediately see that it suffices to prove this existence for abstract complex reflection groups. This is the content of the appendix written by Stephen Griffeth, which concludes the proof of Theorem \ref{fixedpointtheorem G=F}.

\section{Proof of Proposition \ref{prop T0 and PG}}\label{sec tori}

In order to prove Theorem \ref{thm etale quot}, we will need to first prove Proposition \ref{prop T0 and PG}. In this section we show how to deduce it from the results in \cite{ALA}.\\

Lemma \ref{lemma desc Ti} already defines $P_F$ as a complementary $F$-stable subtorus of $T_0$. In order to prove that $P_F$ is an abelian variety, we may reduce to the case where $T=P_F$, (i.e.~$T_0=0$) by restricting the action of $F$ to $P_F$.

The case when the action of $F$ on $T_0(T)$ is irreducible follows from the implication $(1)\Rightarrow(4)$ of \cite[Theorem 1.1]{ALA}. Indeed, although the theorem cited states that the result is for an abelian variety, a look at the proof shows that the only time the existence of an ample line bundle on $T$ was used was in Section 2.1. In this section, if $\sigma\in F$ is a pseudoreflection of order $r$, then the elliptic curve $E_\sigma:=\mathrm{Im}(1-\sigma)$ and abelian subvariety $D_\sigma:=\mathrm{Im}(1+\sigma+\cdots+\sigma^{r-1})$ are defined and shown to be complementary abelian subvarieties with respect to a $G$-invariant polarization. However, the only fact that is actually used in the subsequent proofs is that they are complementary, i.e.~that they generate $T$ and $E_\sigma\cap D_\sigma$ is finite, which is still true for complex tori. In particular, we can conclude that $T$ in this case is isomorphic to a product of elliptic curves, and is thus projective.

As for the reducible case, we resort to \cite[Theorem 2.7]{ALA} where, in the whole section, the projectivity of $T$ is never used when $\dim T^F=0$. We then conclude that this situation splits into products of irreducible cases and use the previous paragraph. This proves Item \eqref{complementary}.

For Item \eqref{fibration}, we use \cite[Proposition 2.9]{ALA}, where once again the projectivity of $T$ can be easily avoided since its proof only uses the complementarity of $T_0$ and $P_F$. We get then the desired fibration and the smoothness of $P_F/F$, which implies by Item \eqref{complementary} that $P_F$ is projective.

\begin{rem}
Given that the intersection $T_0\cap P_F$ measures how far is $T$ from being a direct product (in which case $T/F$ would give the trivial fibration over $T_0$), one could wonder at least what the torsion of $T_0\cap P_F$ looks like. As it can be deduced from the constructions in Section \ref{sec reduction} below, the group $T_0\cap P_F$ can be \emph{any} subgroup of the invariant elements $P_F^F$. The possible groups $P_F^F$ are computed explicitly in the proof of Proposition \ref{prop class admissible irred}. This fact was also used in \cite{Gary} in order to obtain a full classification of these fibrations when the base is an abelian variety.
\end{rem}

\section{Proof of Theorem \ref{thm etale quot}}\label{sec reduction}

Keeping the notations from above, we assume now that $T/F$ is smooth and that $F$ has a fixed point on $T$. In particular, by Proposition \ref{prop T0 and PG}, we have $T=T_0+P_F$ and there is a fibration $\pi:T/F\to T_0/\Delta$ whose fibers are products of projective spaces, where $\Delta$ is the finite group $T_0\cap P_F$. We will briefly recall this construction, which comes from \cite[Prop.~2.9]{ALA}. Following the proof of \cite[Prop.~2.9]{ALA}, we obtain the following commutative diagram, which provides a trivialization the fibration $T/F\to T_0/\Delta$:
\begin{equation}\label{eqn big diagram}
\xymatrix{
T_0\times P_F \ar[d]_F \ar[r]^{\Delta} & T \ar[d]_F \ar[dr]^{G} & \\
T_0\times (P_F/F) \ar[d] \ar[r]^{\Delta} & T/F \ar[d] \ar[r]_{G/F} & T/G \\
T_0 \ar[r]^{\Delta} & T_0/\Delta. }
\end{equation}
Here, a group besides an arrow denotes the quotient map by an action of that group. Moreover, again by Proposition \ref{prop T0 and PG}, we know that $P_F/F\simeq \mathbb{P}^{n_1}\times\cdots\times\mathbb{P}^{n_r}$.

Let us prove now that $G/F$ acts naturally on $T_0/\Delta$. First of all, we show that $G_\lin$ preserves $T_0$ and $P_F$, which amounts to proving that it preserves the corresponding tangent spaces. Now, these tangent spaces correspond to the decomposition of $T_0(T)$ into trivial and non trivial representations of $F_\lin$ respectively. Since $F_\lin$ is normal in $G_\lin$ it is an easy exercise to check that $G_\lin$ preserves these two subspaces.

Having proved this, we note that the arrow $p:T\to T_0/\Delta$ actually corresponds to the quotient by the subtorus $P_F$. Then, since $G_\lin$ preserves $P_F$, we see that its action descends naturally to an action on $T_0/\Delta$, which is easily seen to be trivial on $F_\lin$, hence it corresponds to an action of $G_\lin/F_\lin$. Define then, for $x+\Delta\in T_0/\Delta$ and $gF\in G/F$,
\[(gF)(x+\Delta):=g_\lin(x)+p(g(0))+\Delta.\]
A direct computation (using the fact that $G$ acts on $T$) proves that this defines an action of $G/F$. Moreover, by its very definition, this action is $G$-equivariant with respect to $p$, hence $(G/F)$-equivariant with respect to the arrow $T/F\to T_0/\Delta$, which factors $p$ in the diagram above.

Defining then $M:=(T_0/\Delta)/(G/F)$, we can complete the lower-right-side of the diagram into a commutative square
\begin{equation}\label{diagram}\xymatrix{
T/F \ar[r]_{G/F} \ar[d] & T/G \ar[d] \\
T_0/\Delta \ar[r]_{G/F} & M.}
\end{equation}
Since we know that the action above is free, all we are left to prove is that the action of $G/F$ on $T_0/\Delta$ is free as well. This immediately implies that the arrow on the right is a fibration having the desired properties.

Consider then $b\in T_0/\Delta$ and let $S:=\stab_{G/F}(b)\leq G/F$. We must prove that $S$ is the trivial group. If we denote $X:=(T/F)$ and $X_b$ the fiber over $b$, then $S$ acts naturally on $X_b$, which is isomorphic to a product of projective spaces. Now we need the following proposition, which follows from \cite[Thm.~1]{alvaro}:

\begin{proposition}\label{prop aut prod Pn}
For any positive integers $m_1,n_1,\ldots,m_r,n_r$ with the $n_i$ pairwise different, we have 
\[\aut\left(\prod_{i=1}^r(\mathbb{P}^{n_i})^{m_i}\right)\simeq\prod_{i=1}^r(\mathrm{PGL}_{n_i+1}(\C)^{m_i}\rtimes S_{m_i}).\]
\end{proposition}

From this, it is easy to see that every element in $\aut(X_b)$ must fix a point in $X_b$, which implies that every element in $S$ fixes a point in $T/F$. By Item (1) of Theorem \ref{fixedpointtheorem}, we see that such an element must be the identity, which concludes the proof.

\section{Reducing the classification of smooth quotients to the hyperelliptic case}\label{sec red hyperell case}
In this section, we provide a method that allows to classify smooth quotients of complex tori (resp.~abelian varieties) fibred over a given torus (resp.~abelian variety) or hyperelliptic manifold (resp.~variety). More precisely, we pin down the automorphisms of a product of projective spaces that can appear as glueing morphisms in such a fibration.\\

First of all, let $A$ be a complex torus, let $G\leq\Lin(A)$ be a finite group of automorphisms such that $A/G\simeq\mathbb{P}^n$ (in particular, $A$ is an abelian variety), and let $\pi:A\to \mathbb{P}^n$ be the quotient map. Note that a morphism $\varphi\in\aut(A)$ descends to an automorphism of $\mathbb{P}^n$ if and only if it belongs to the normalizer of $G$ in $\aut(A)$. In particular, there is a natural group homomorphism
\[\pi_*:N_{\aut(A)}(G)\to\aut(\mathbb{P}^n)=\mathrm{PGL}_{n+1}(\C).\]

\begin{definition}
We say that a group $\Gamma$ of automorphisms of $\bb P^n$ is \emph{admissible} if there exists a pair $(A,G)$ where $A$ is an abelian variety and $G$ is a finite subgroup of $\Lin(A)$ such that $A/G\simeq\bb P^n$ and $\Gamma$ is in the image of $\pi_*$.
\end{definition}

Given that pairs $(A,G)$ as above were classified in \cite{ALA}, we get a classification of admissible groups of automorphisms.

\begin{proposition}\label{prop class admissible irred}
Let $\Gamma\leq\aut(\bb P^n)\cong\mathrm{PGL}_{n+1}(\C)$ be admissible. Then, up to conjugacy, $\Gamma$ is contained in one of the following subgroups:
\begin{itemize}
\item[$(A_n)$] The group $(\Z/(n+1)\Z)^2\rtimes (\Z/m\Z)$, with $m=4,6$, which satisfies the following conditions:
\begin{itemize}
\item the generators $\alpha_1,\alpha_2\in (\Z/(n+1)\Z)^2$ fix exactly $n+1$ points each;
\item $\mathrm{Fix}(\alpha_1)\cap \mathrm{Fix}(\alpha_2)=\emptyset$;
\item the generator $\alpha_3$ of $\Z/m\Z$ satisfies:
\begin{itemize}
\item $\alpha_3\alpha_1\alpha_3^{-1}=\alpha_2$ and $\alpha_3\alpha_2\alpha_3^{-1}=\alpha_1^{-1}$ if $m=4$;
\item $\alpha_3\alpha_1\alpha_3^{-1}=\alpha_2$ and $\alpha_3\alpha_2\alpha_3^{-1}=\alpha_1^{-1}\alpha_2$ if $m=6$.
\end{itemize}
\end{itemize}
\item[$(B_n)$] The finite groups $F=D_3,\,D_4,\,A_4$, seen as subgroups of $\mathrm{PGL}_{n+1}(\C)$ in the following way: take \underline{the} representation of $F$ in $\mathrm{PGL}_2(\C)$, make it act on $(\bb P^1)^n$ diagonally and consider the induced action on the symmetric product $\mathrm{Sym}^n(\bb P^1)=\bb P^n$.
\end{itemize}
\end{proposition}

We see then that there are essentially 5 different types of maximal admissible groups of automorphisms.

\begin{rem}
The group described in case $(A_n)$ is indeed unique up to conjugacy. The condition of having exactly $n+1$ fixed points amounts to the fact that a representative matrix in $\mathrm{GL}_n$ has every $(n+1)$-st root of unity as an eigenvalue. Using moreover the fact that $\alpha_1$ and $\alpha_2$ commute and that their fixed point sets do not intersect, one easily gets that, up to conjugacy:
\begin{align*}
\alpha_1:[x_0:x_1:\cdots : x_n]&\mapsto [x_0:\zeta_{n+1} x_1:\cdots :\zeta_{n+1}^nx_n],\\
\alpha_2:[x_0:x_1:\cdots :x_n] &\mapsto [x_n:x_0:\cdots : x_{n-1}].
\end{align*}
Once these two elements are fixed, it takes a tiresome but direct computation to see that elements satisfying the property of $\alpha_3$ are unique up to multiplication by elements in $\langle\alpha_1,\alpha_2\rangle$. For instance, the element $\alpha_3$ for $n=4$ maps $[x_0:x_1:\cdots :x_n]$ to $[y_0:y_1:\cdots:y_n]$ with $y_i:=\sum_{j=0}^n\zeta_{n+1}^{-ij}x_j$.

In case $(B_n)$ it is also direct to obtain explicit generators. Generators as subgroups of $\pgl_2(\C)$ are given by:
\[D_m=\left\langle\begin{pmatrix} 0 & 1 \\ 1 & 0\end{pmatrix},\begin{pmatrix}\zeta_m & 0 \\ 0 & 1\end{pmatrix}\right\rangle,\qquad A_4=\left\langle\begin{pmatrix} i & 0 \\ 0 & -i\end{pmatrix},\begin{pmatrix} 0 & 1 \\ -1 & 0\end{pmatrix},\begin{pmatrix} \frac{1+i}{2} & \frac{-1+i}{2} \\ \frac{1+i}{2} & \frac{1-i}{2}\end{pmatrix}\right\rangle\]
Then, a direct computation gives that the generators of the action of, say, $D_3$ on $\bb P^n$ are:
\begin{align*}
\beta_1:[x_0:x_1:\cdots : x_n]&\mapsto [x_n:x_{n-1}:\cdots :x_0],\\
\beta_2:[x_0:x_1:\cdots :x_n] &\mapsto [x_0:\zeta_3x_1:\cdots :\zeta_3^nx_n].
\end{align*}
\end{rem}

\begin{proof}[Proof of Proposition \ref{prop class admissible irred}]
Consider the pair $(A,G)$ such that elements in $\Gamma$ lift to automorphisms of $A$ via the quotient morphism. Since the preimage of $\Gamma$ is contained in the normalizer $N_{\aut(A)}(G)$ of $G$ in $\aut(A)$, we only need to compute these normalizers for every pair $(A,G)$ such that $A/G\simeq\bb P^n$ and study their images in $\aut(\bb P^n)=\mathrm{PGL}_{n+1}(\C)$. We claim that these correspond to the subgroups on the list above.

First of all, let $\phi\in\aut(A)$ be defined as $\phi(x)=\phi_\lin(x)+t$ with $t\in A$ and $\phi_\lin\in\Lin(A)$ and assume that $\phi$ normalizes $G$. We compute then, for $x\in A$ and $g\in G$,
\begin{align*}
(\phi^{-1}g\phi)(x) &=(\phi^{-1}g)(\phi_\lin(x)+t)\\
&=\phi^{-1}((g\phi_\lin)(x)+g(t))\\
&=\phi_\lin^{-1}((g\phi_\lin)(x)+g(t)-t)\\
&=(\phi_\lin^{-1}g\phi_\lin)(x)+\phi_\lin^{-1}(g(t)-t).
\end{align*}
Applying this to $x=0$ we see that $g(t)-t=0$ since $\phi^{-1}g\phi\in G\subseteq\Lin(A)$. Hence $t\in A^G$. Taking then $t\in A^G$, we see that $\phi_\lin$ must normalize $G$ as well. In other words,
\[N_{\aut(A)}(G)\subseteq A^G\rtimes N_{\Lin(A)}(G)\subseteq A\rtimes\Lin(A)=\aut(A).\]
On the other hand, the computations above confirm that the first inclusion is in fact an equality. Thus, we are only left to compute $A^G$ and the quotient $N/G$, where $N$ is the normalizer of $G$ in $\Lin(A)$, then compute the action of $A^G\rtimes (N/G)$ on the quotient $A/G\simeq\bb P^n$. Note by the way that this action is faithful since the action of $A^G\rtimes N$ on $A$ is.

By \cite[Thm.~1.1]{ALA}, we know that there are essentially two families of pairs $(A,G)$, which will give us the subgroups of type $(A_n)$ and $(B_n)$, respectively. There is also an exceptional case in dimension 2 by \cite[Thm.~1.1]{ALAQ}, for which we will prove that it does not give rise to new subgroups. We go then case by case.\\

{\bf The $(A_n)$ case:} Assume that $G=S_{n+1}$ with $n\geq 2$ acting on
\[A=\{(x_1,\ldots,x_{n+1})\in E^{n+1}\mid \sum x_i=0\},\]
in the obvious way. It is easy to see then that
\[A^G=\{(x,\ldots,x)\in E^{n+1}\mid x\in E[n+1]\}\simeq(\Z/(n+1)\Z)^2.\]
This gives the generators $\alpha_1$ and $\alpha_2$. Indeed, one can easily prove that each generator of $E[n+1]$ fixes exactly $n+1$ points in $\bb P^n$ and that these differ for different generators.

On the other hand, it is well-known that the normalizer of $S_{n+1}$ in $\mathrm{GL}_{n+1}(\C)$ (seen as permutation matrices) corresponds to the product of $S_{n+1}$ itself with its centralizer, which is the subgroup given by the matrices of the form
\begin{equation}\label{eq cent Sn}
M=\begin{pmatrix} a & b & \cdots & b \\ b & a & & \vdots \\ \vdots & & \ddots & b \\ b & \cdots & b & a\end{pmatrix},\text{ with }a,b\in\C\text{ and }\det(M)\neq 0.
\end{equation}
If we restrict the action of such a matrix to $T_0(A)\subseteq T_0(E^{n+1})=\C^{n+1}$, which corresponds to the hyperplane of trivial sum, we see that it induces multiplication by $a-b$. Since this multiplication must be invertible on $A\simeq E^n$, we are only left with roots of unity, assuming $E$ admits such a multiplication. Such elements correspond to the group $N/G$ and give the generator $\alpha_3$ when we consider multiplication by $\zeta_4$ and $\zeta_6$ respectively (of course, multiplications by $\zeta_3$ and $-1$ are contained in these).\\

{\bf The $(B_n)$ case:} Assume now that $G=G(m,1,n)=C^n\rtimes S_n$ with $C$ a subgroup of $\Lin(E)$ of order $m\in\{2,3,4,6\}$. It is easy to see then that
\[A^G=\{(x,\ldots,x)\in E^n\mid x\in E^C\}\simeq E^C.\]
On the other hand, note that $\Lin(A)=\mathrm{GL}_n(\Z[\zeta_\ell])$ with $\ell$ the order of $\Lin(E)$. Here $S_n$ corresponds to permutation matrices and $C^n$ to diagonal matrices. Consider then an element $h\in N$. Since $S_n$ is a subgroup of $G$, we have that $hS_n h^{-1}$ must correspond to a certain section of $S_n$ in $G$. Now sections of $S_n$ in $G=C^n\rtimes S_n$ are classified, up to conjugation in $G$, by the cohomology group $H^1(S_n,C^n)$ (cf.~\cite[I.2, Exer.~1]{NSW}). And since $C^n$ is clearly an induced $S_n$-module (it admits a basis permuted by the group $S_n$), we know by Shapiro's Lemma (cf.~\cite[Prop.~1.6.4]{NSW}) that $H^1(S_n,C^n)\simeq H^1(S_{n-1},C)$ with $S_{n-1}$ acting trivially on $C$. Then we have
\[H^1(S_{n-1},C)=\hom(S_{n-1},C)=\begin{cases} \Z/2\Z &\text{if } m\text{ is even,} \\ \{0\} & \text{if } m \text{ is odd.}\end{cases}\]
So there is at most one class of sections that is not conjugate to the original subgroup $S_n$. A section of $S_n$ corresponding to this nontrivial class is given by elements of the form $(-I)^{\mathrm{sgn(\sigma)}}\sigma$, where $I$ is the identity matrix, $\sigma$ is a permutation matrix and $\mathrm{sgn}(\sigma)$ is its sign. Indeed, these are contained in $G$ for even $m$ (since $-I$ is) and correspond to a complex representation of $S_n$ that is clearly different from the original one (actually, it is the tensor product with the sign representation). We see then that such a subgroup cannot be obtained by conjugation by \emph{any} matrix in $\mathrm{GL}_n(\C)$ and so $hS_nh^{-1}$ corresponds to a section conjugate in $G$ to $S_n$. Thus, up to modifying our element $h$ by an element in $G$, we may assume that it centralizes $S_n$ and it is thus given by \eqref{eq cent Sn} if $n>2$. If $n=2$, then the centralizer in $\mathrm{GL}_2(\C)$ is slightly bigger and is given by matrices of the form
\[M=\begin{pmatrix} a & b \\ \pm b & \pm a\end{pmatrix},\text{ with }a,b\in\C\text{ and }\det(M)\neq 0.\]
Recalling then that $h\in\Lin(A)=\mathrm{GL}_n(\Z[\zeta_\ell])$, we see that the determinant must be a unit in $\Z[\zeta_\ell]$. This gives $b=0$ for $n>2$ and finitely many matrices to check for $n=2$. After a quick check and possibly modifying $h$ by an element in $G$, we get that $h$ corresponds to multiplication by a root of unity, assuming $E$ admits such a multiplication.

With all this, and recalling that the quotient $A/G$ factors through the quotient $A/(C^n)=(E/C)^n\simeq(\bb P^1)^n$, we see that:
\begin{itemize}
\item If $m=2$, we have $E^C=E[2]\simeq(\Z/2\Z)^2$ and, if $E$ has multiplication by either $\zeta_4$ or $\zeta_6$, then multiplication by this element gives an automorphism in $N$ that is not in $G$, so that $N/G$ is either $\Z/2\Z$ or $\Z/3\Z$. Then $A^G\rtimes(N/G)$ is either $D_4$ or $A_4$ acting faithfully and diagonally on the product $A/(C^n)\simeq (\bb P^1)^n$.
\item If $m=3$, we have $E^C\simeq\Z/3\Z$ and $E$ has multiplication by $\zeta_6$, which gives an automorphism in $N$ that is not in $G$, so that $N/G\simeq \Z/2\Z$. Then $A^G\rtimes(N/G)$ corresponds to $S_3$ acting faithfully and diagonally on the product $A/(C^n)\simeq (\bb P^1)^n$.
\item If $m=4$, we have $E^C\simeq\Z/2\Z$ and $E$ has multiplication by $\zeta_4$, but such a multiplication is already contained in $G$, so that $A^G\rtimes(N/G)$ is simply $\Z/2\Z$ and it can be seen as contained in one of the groups above.
\item If $m=6$, we have $E^C=\{0\}$ and $E$ has multiplication by $\zeta_6$, but such a multiplication is already contained in $G$, so that $A^G\rtimes(N/G)$ is trivial.
\end{itemize}
These actions commute with the action of $S_n$ on $(\bb P^1)^n$ and thus give the desired action on $\bb P^n$. Note by the way that this last case-by-case analysis also works for $n=1$.\\

{\bf The exceptional case:} Here we have $n=2$, $A=E^2$ where $E=\C/\Z[i]$ and $G$ is the group of order 16 generated by the matrices
\[\left\{\begin{pmatrix} -1 & 1+i \\ 0 & 1\end{pmatrix}\right.,\, \begin{pmatrix} -i & i-1 \\ 0 & i\end{pmatrix},\, \left.\begin{pmatrix} -1 & 0 \\ i-1 & 1\end{pmatrix} \right\}\]
that act on $A$ in the obvious way. A direct computation gives that $\{\pm I,\pm i I\}\subset G$, from which one immediately obtains
\[A^G=\langle({\textstyle \frac{1+i}{2},0),(0,\frac{1+i}{2}})\rangle\simeq (\Z/2\Z)^2.\]
As for the normalizer $N$ of $G$ in $\Lin(A)$, recall once again that $\Lin(A)=\mathrm{GL}_2(\Z[i])$. We therefore have a natural homomorphism 
\[\pi:\Lin(A)\to\mathrm{PGL}_2(\C)\]
whose kernel consists of the matrices $\{\pm I,\pm i I\}\subset G$, and whose image $\pi(G)\subseteq\mathrm{PGL}_2(\C)$ is isomorphic to the Klein 4-group. The normalizer $N$ of $G$ is then sent to a finite subgroup of the normalizer of $\pi(G)$ in $\mathrm{PGL}_2(\C)$.

It is well-known that two isomorphic finite subgroups of $\mathrm{PGL}_2(\C)$ belong to the same conjugacy class. Moreover, the finite subgroups of $\pgl_2(\C)$ are either cyclic, dihedral, $A_4$, $S_4$ or $A_5$. Out of these, the only ones having a normal subgroup isomorphic to the Klein 4-group are the Klein group itself, $A_4$ and $S_4$, thus the image of $N$ in $\pgl_2(\C)$ must be contained in $S_4$.

We observe that the matrix 
\[M:=\begin{pmatrix}i & 1\\0 & 1\end{pmatrix}\in\Lin(A)\] 
is in $N$ since the subgroup $H:=\langle M,G\rangle\leq\Lin(A)$ is of order $32$. In particular, the group $\pi(H)$ is of order 8, and so $\pi(N)$ is either of order $8$ or $24$. On the other hand, a simple calculation shows that the subgroup of $\pgl_2(\C)$ generated by $\pi(G)$ and the matrix
\[S:=\begin{pmatrix} -2i-2 & i-1 \\ -2 & 2i\end{pmatrix}\]
is \underline{the} subgroup of $\pgl_2(\C)$ isomorphic to $S_4$ that contains $\pi(G)$, and since $S\notin\pi(\Lin(A))$, we conclude that $\pi(N)=\pi(H)$ and thus $N=H$.

Putting everything together, we see that the group $A^G\rtimes(N/G)$ is isomorphic to $D_4$ and acts faithfully on $\bb P^2$. But $D_4$ has only one faithful representation in $\pgl_3(\C)$ (this is an easy exercise on linear representations of $D_4$), so it is already considered in case $(B_2)$.
\end{proof}

We will now extend the notion of an admissible group action to the more general setting of products of projetive spaces.

\begin{definition}\label{defi admissible general}
Let $\Gamma$ be a finite group acting faithfully as a group of automorphisms of a product of projective spaces
\[X:=\mathbb{P}^{n_1}\times\cdots\times\mathbb{P}^{n_r}.\]
We say that $\Gamma$ is \emph{admissible} if there exists a pair $(A,G)$ with
\[A=A_1\times\cdots\times A_r\quad\text{and}\quad G=G_1\times\cdots\times G_r,\]
where $A_i$ is an abelian variety, $G_i$ is a subgroup of $\Lin(A_i)$ such that $A_i/G_i\simeq\mathbb{P}^{n_i}$, and each automorphism in $\Gamma$ lifts to an automorphism of $A$ via the projection $\pi:A\to X=A/G$.\\
Moreover, we define $\pi^*\Gamma:=(\pi_*)^{-1}(\Gamma)$, where $\pi_*$ is the group homomorphism
\[\pi_*:N_{\aut(A)}(G)\to\aut(X).\]
\end{definition}

We note that this case is more complicated than the irreducible case (i.e.~when $r=1$), since we can also have permutations of factors. Fortunately, Proposition \ref{prop aut prod Pn} tells us that this in fact the only possible new addition to admissible automorphisms. We obtain thus a full classification of admissible automorphisms of products of projective spaces.

\begin{proposition}\label{prop class admissible general}
Let $\Gamma$ be a finite group acting faithfully on
\[X:=\mathbb{P}^{n_1}\times\cdots\times\mathbb{P}^{n_r}.\]
Then $\Gamma$ is admissible if and only if there is a partition $\cal P=\{P_1,\ldots,P_s\}$ of the set $\{1,\ldots,r\}$ such that $\Gamma$ is a subgroup of
\[\prod_{P\in \cal P}\left(\Big(\prod_{i\in P} \Gamma_i\Big)\rtimes S_{P}\right)=\left(\prod_{i=1}^r\Gamma_i\right)\rtimes\left(\prod_{P\in\cal P} S_P\right),\]
where the $\Gamma_i$ are maximal admissible subgroups of $\pgl_{n_i+1}(\C)$ and, for $i,j\in P_k$, $\Gamma_i$ and $\Gamma_j$ have \emph{the same type} (in particular, $n_i=n_j$) and the symmetric groups $S_P$ permute the variables in $P$.
\end{proposition}

\begin{proof}
Note first that a group as the one in the statement of the proposition is admissible. Indeed, since for a fixed $P\in\cal P$ the group the $S_P$ permutes groups $\Gamma_i$ of the same type, we may assume that the pairs $(A_i,G_i)$ and $(A_j,G_j)$ are isomorphic for $i,j\in P$. Then clearly the elements in $S_P$ lift to automorphisms of $\prod_{i\in P} A_i$ and, since each $\Gamma_i$ is admissible, we get a lift of the whole group $(\prod_{i\in P} \Gamma_i)\rtimes S_{P}$ to $\aut(\prod_{i\in P} A_i)$. Since this is true for every $P$, we get the claim.\\

On the other hand, by Proposition \ref{prop aut prod Pn} we know that
\[\aut(X)=\left(\prod_{i=1}^r\pgl_{n_i+1}(\C)\right)\rtimes S,\]
where $S$ is a product of symmetric groups permuting the components of same dimension. In particular, every element in $\aut(X)$ can be written as a composition $\psi\sigma$ with $\psi\in \prod_{i=1}^r\pgl_{n_i+1}(\C)$ and $\sigma$ a permutation in $S\leq S_r$.

It is clear then from Proposition \ref{prop class admissible irred} that the liftable elements of the product $\prod_{i=1}^r\pgl_{n_i+1}(\C)$ must be contained in a product of maximal admissible subgroups, hence a product of $\Gamma_i$'s as in the statement of the proposition.

We claim now that, in order to lift an automorphism sending the $i$-th component to the $j$-th component (i.e.~of the form $\psi\sigma$ with $\sigma$ sending $i$ to $j$), the pair $(A_i,G_i)$ must be isomorphic to $(A_j,G_j)$. Assuming the claim, the corresponding maximal admissible subgroups $\Gamma_i$ and $\Gamma_j$ must be of the same type since different types arise from non isomorphic pairs $(A,G)$, as can be deduced from the proof of Proposition \ref{prop class admissible irred}. These isomorphisms can then be used to construct liftable automorphisms in $\aut(X)$ of the form $\psi'\sigma$ for some $\psi'\in \prod_{i=1}^r\pgl_{n_i+1}(\C)$ and $\sigma\in S$. Then the whole coset $(\prod_{i=1}^r\Gamma_i)\psi'\sigma$ is contained in an admissible subgroup. Thus, our starting $\psi$ must belong to $(\prod_{i=1}^r\Gamma_i)\psi'$ since otherwise we contradict the maximality of the $\Gamma_i$'s. This implies that $\psi{\psi'}^{-1}$ is liftable and we can rewrite our starting automorphism as $(\psi{\psi'}^{-1})(\psi'\sigma)$, which is clearly the unique way of writing it as a product of liftable automorphisms. This implies that our admissible group $\Gamma$ is contained in a semi-direct product like the one in the statement of the proposition, so we are reduced to proving the claim.\\

Consider then a lift $\varphi\in\aut(A)$, with $A=\prod_{i=1}^rA_i$, of an automorphism $\psi\sigma\in\aut(X)$, where $\sigma$ sends $i$ to $j$. Let $\pi_i:A_i\to\mathbb{P}^{n_i}$ and $\pi:=\pi_1\times\cdots\times\pi_r:A\to X$ be the natural projections. Let $(a_1,\ldots,a_r)\in X$ be an $S$-invariant point and, for each $i$, define the inclusions
\[\kappa_{i}:\mathbb{P}^{n_i}\to X\,:\,x\mapsto (a_1,\ldots,a_{i-1},x,a_{i+1},\ldots,a_r).\]
Write $\psi=(\psi_1,\ldots,\psi_r)$ and let $b_i:=\psi_i(a_i)$. Then \[\psi\sigma\kappa_i(\bb P^{n_i})=\{(b_1,\ldots,b_{j-1},x,b_{j+1},\ldots,b)\mid x\in\bb P^{n_j}\}.\]
Thus, $\pi^{-1}(\psi\sigma\kappa_i(\bb P^{n_i}))$ is the disjoint union of sets of the form
\begin{equation}\label{eq A_j}
\{(z_1,\ldots,z_{j-1},x,z_{j+1},\ldots,z_r)\mid x\in A_j\}\simeq A_j,
\end{equation}
where the $z_i$'s run through the corresponding preimages of the $b_i$'s via $\pi_i$. By the same argument, the preimage of $\kappa_i(\bb P^{n_i})$ via $\pi$ is the disjoint union of subvarieties isomorphic to $A_i$, and these are sent by $\varphi$ to each of the components in \eqref{eq A_j} by connectedness and thus $A_i\simeq A_j$. Since moreover the projections from these components to $\bb P^{n_i}$ and $\bb P^{n_j}$ correspond to quotients by the respective actions of $G_i$ and $G_j$ on $A_i$ and $A_j$, we see that we must have $G_i\simeq G_j$. But this means that $(A_i,G_i)\simeq (A_j,G_j)$ since each pair is uniquely determined by the dimension of the abelian variety and the order of the group. 
\end{proof}

With this classification at hand, we may proceed to a classification of smooth quotients of complex tori fibred over a given base (which will be an \'etale quotient of a complex torus). This can be seen as a sort of bridge between the case of quotients of complex tori by finite groups that fix the origin and hyperelliptic manifolds.

\begin{theorem}\label{thm classification}
Let $Z_0$ be a complex torus and write $Z=Z_0/H$ with $H$ a finite group acting freely (possibly with translations). Consider further a product of projective spaces
\[X:=\mathbb{P}^{n_1}\times\cdots\times\mathbb{P}^{n_r},\]
and let $H$ act on $X$ via admissible automorphisms, i.e.~we fix a morphism $H\to\aut(X)$ with admissible image. Define $Y:=(Z_0\times X)/H$, which comes with a natural arrow $Y\to Z=Z_0/H$. Then the arrow $Y\to Z$ is a fibration of products of projective spaces and the manifold $Y$ corresponds to a smooth quotient of a complex torus.

Moreover, every smooth quotient of a complex torus that fibers on products of projective spaces over $Z$ as in Theorem \ref{thm etale quot} can be obtained with this construction.
\end{theorem}

\begin{proof}
Note that $Y$ is smooth since the action of $H$ on $Z_0$ is free, and hence so is the case for $Z_0\times X$. It is evident by construction that the fibers of the arrow $Y\to Z$ are isomorphic to $X$, hence it is indeed a fibration of products of projective spaces. Since $H$ acts via an admissible group on $X$, there exists a pair $(A,F)$ with $A$ an abelian variety, $F$ a finite group and $A/F=X$. Recalling Definition \ref{defi admissible general}, we can consider $G:=\pi^*H\subset\aut(A)$. We make $G$ act on the product $Z_0\times A$ as follows: the action on $Z_0$ is given by the arrow $\pi:G\to H$ and the action of $H$ on $Z_0$, whereas the action on $A$ is the obvious one since $G\subset\aut(A)$. Since the kernel of $\pi:G\to H$ is clearly $F$, which acts trivially on $Z_0$, we see that
\[(Z_0\times A)/G=((Z_0\times A)/F)/H=(Z_0\times(A/F))/H=(Z_0\times X)/H=Y.\]
This proves the first assertion of the Theorem.\\

Assume now that there is a fibration $Y\to Z$ of products of projective spaces arising from a smooth quotient of a complex torus as in Theorem \ref{thm etale quot}, i.e. $Y=T/G$ for some complex torus $T$ and some finite group $G$ and the fibers of $Y\to Z$ are isomorphic to $X$. We may consider then the following diagram, which follows from diagrams \eqref{diagram} and \eqref{eqn big diagram} from Section \ref{sec reduction}:
\[\xymatrix{
T_0\times P_F \ar[d]_F \ar[r]^{\Delta} & T \ar[d]_F \ar[dr]^{G} & \\
T_0\times X \ar[d] \ar[r]^{\Delta} & T/F \ar[d] \ar[r]_{G/F} & Y \ar[d] \\
T_0 \ar[r]^{\Delta} & T_0/\Delta \ar[r]_{G/F} & Z.}
\]
Here we have $X=P_F/F$ and the diagram is giving us a trivialization of the fibration $Y\to Z$ as $T_0\times X\to T_0$. Note that both arrows $T_0\to T_0/\Delta$ and $T_0/\Delta\to Z$ are \'etale Galois morphisms. We may consider then the Galois closure $Z_0\to Z$ of the composition $T_0\to Z$, which is also an \'etale cover of $T_0$ and thus $Z_0$ is a complex torus. Denoting by $H$ the Galois group of $Z_0\to Z$ and pulling back the trivialization from $T_0$ to $Z_0$, we see that $Y\to Z$ is obtained as the quotient of $Z_0\times X\to Z_0$ by $H$ for some action of $H$ on $X$. Using the diagram above, it is easy to see that this action is admissible, which proves the second assertion and concludes the proof of the theorem.
\end{proof}

\begin{rem}
Note that already passing from $T_0$ to $T_0/\Delta$ we are considering translations. And it is clear that the morphism $Z_0\to T_0$ consists of translations as well. This is why we need to present the quotient $Z=Z_0/H$ as a quotient possibly including translations. In fact, already in \cite[Prop.~2.9]{ALA} we get non trivial fibrations over abelian varieties, which are trivialized by considering isogenies, i.e.~actions via translations.

Note however that, if one wishes to actually compute \emph{all} possible fibrations for a given base manifold $Z$ and fiber manifold $X$, a trivial bound for the degree of the isogenies that must be considered is the order of the maximal admissible subgroups of $\aut(X)$, since the action must factor through there.
\end{rem}

\section{Examples of non-abelian hyperelliptic manifolds}\label{sec examples}
In this section we present a family of pairs $(T,G)$ such that $T/G$ is an \'etale quotient and $G$ is non-abelian and arbitrarily large. These are intended to show that the situation is much richer (and hence more complicated) than the already treated case of actions fixing the origin. Recall that we assume that $G$ contains no translations.\\

Before we start, we make the following easy remarks:

\begin{enumerate}
\item A finite \emph{abelian} group $G$ can always act freely on some complex torus $T$. Indeed, since the $n$-torsion of a complex torus of dimension $d$ is isomorphic to $(\Z/n\Z)^{2d}$, we can always embed $G$ into $T$ granted that the dimension of the latter is big enough. This induces a natural action of $G$ on $T$ by translations, hence free.
\item A finite group $G$ can always act faithfully on some complex torus $T$ fixing the origin. Indeed, it suffices to take $|G|$ copies of any complex torus $Z$ and let $G$ act on them by permutation.
\item Given a non trivial action of $G$ on $T$ fixing the origin, there exist finite sub-$G$-modules of $T$ with non trivial action. Indeed, one can simply take the $n$-torsion subgroup $T[n]$ for $n$ large enough.
\end{enumerate}

Given the first two obvious facts stated above, there is an easy way of generating free actions of a finite \emph{abelian} $G$ on some complex torus containing no translations: just consider the direct product of the two situations. The action by translations ensures the freeness and the faithful action fixing the origin ensures that there are no translations for the action on the product. Note that this is precisely the way in which bielliptic surfaces are constructed.

One could wonder then whether there are pairs $(T,G)$ with \'etale quotient and \emph{non-abelian} $G$. This can be done via the following construction, which also uses the trivial remarks from above.

\begin{itemize}
\item[(1)] Fix a complex torus $T_0$ and a finite subgroup $G_1\subseteq T_0$. Then $G_1$ acts freely on $T_0$ by translations.
\item[(2)] Fix a complex torus $Z_1$ with a faithful action of $G_1$ that fixes the origin. Then $G_1$ acts freely on $T_1:=T_0\times Z_1$ with no translations.
\item[(3)] Fix now a finite sub-$G_1$-module $F_1\subseteq Z_1$ with non trivial action and define $G_2:=F_1\rtimes G_1$. Then $G_2$ acts on $T_1=T_0\times Z_1$ naturally: $F_1$ acts by translations and $G_1$ as in the previous step. One can check that the action is well-defined by a simple computation. Moreover, the action is easily seen to be free, but it contains translations.
\item[(4)] Fix a complex torus $Z_2$ with a faithful action of $G_2$ that fixes the origin. Then $G_2$ acts freely on $T_2:=T_1\times Z_2$ with no translations.
\end{itemize}

The pair $(T_2,G_2)$ is then an example with $T_2/G_2$ an \'etale quotient and non-abelian $G_2$. If one wishes to go further, there is always:

\begin{itemize}
\item[(5)] Iterate steps (3) and (4) at will replacing $T_i$, $Z_i$, $F_i$, and $G_i$ by $T_{i+1}$, $Z_{i+1}$, $F_{i+1}$, $G_{i+1}$ as needed.
\end{itemize}

This recipe tells us that we can construct arbitrarily large iterated semi-direct products of abelian groups with free action on some complex torus with no translations. In particular, such a pair $(T,G)$ would give an \'etale quotient $T/G$.

Note that the analytic representation of $G_\lin$ induced by all of these examples always has a copy of the trivial representation as a direct factor. Much more interesting examples in which there is no such copy can be found in \cite{CatDem} and \cite{Bruno}.

\begin{appendix}

\section{A result on complex reflection groups, by Stephen Griffeth}

The purpose of this appendix is to prove the result mentioned in the proof of Proposition \ref{prop fixed point} above: each complex reflection group contains an element not contained in any proper reflection subgroup (Theorem \ref{appendixthm} below). We first fix notation and definitions, and then give the proof. One may assume the group is irreducible and hence use the classification of irreducible reflection groups. The proof is then essentially uniform for the well-generated groups but requires a case-by-case check for the remaining groups. It would of course be interesting to find a completely uniform proof.

\subsection{Reflection groups: Coxeter numbers}

 Let $\hh$ be a finite-dimensional $\C$-vector space and let $W \subseteq \mathrm{GL}(\hh)$ be a finite group of linear transformations of $\hh$. The set of \emph{reflections} in $W$ is
$$R=\{r \in W \ | \ \mathrm{codim}(\mathrm{fix}(r))=1 \},$$ and $W$ is a \emph{reflection group} if it is generated by $R$. We write $$\mathcal{A}=\{\mathrm{fix}(r) \ | \ r \in R \}$$ for the set of reflecting hyperplanes for the reflections in $W$. 

For a linear group $W$, the ring $\C[\hh]^W$ of $W$-invariant polynomial functions on $\hh$ is a polynomial ring if and only if $W$ is a reflection group. In this case, the sequence of \emph{degrees} of $W$ is defined by $d_1\leq d_2 \leq \cdots \leq d_n$ where $n=\mathrm{dim}(\hh)$ and $d_1,\dots,d_n$ are the degrees of homogeneous polynomials generating $\C[\hh]$. 

When $\hh$ is an irreducible $\C W$-module we call $W$ \emph{irreducible}. For an irreducible reflection group $W$, the element 
$$\sum_{r \in R} (1-r)$$ is central in $\C W$ and hence acts as a scalar on $\hh$; we call this scalar the \emph{Coxeter number} $h$ of $W$. Since it is the value of a central character on an integer linear combination of class sums, it is an integer, and direct calculation with the trace shows
$$h=\frac{N+N^*}{n},$$ where $N=|R|$ is the number of reflections in $W$, $N^*=|\mathcal{A}|$ is the number of reflecting hyperplanes for $W$, and $n=\mathrm{dim}(\hh)$ is the rank of $W$. We remark that this number $h$ arises also (and for the same reason) in the context of the representation theory of rational Cherednik algebras and Catalan numbers for complex reflection groups; see e.g. \cite{Gri} and \cite{GoGr}. An irreducible complex reflection group $W$ is called \emph{well-generated} if there is a subset of $R$ of cardinality precisely the rank $n$ of $W$ that generates $W$. 

Since $W$ is a finite group, there exists a positive definite Hermitian form on $\hh$ which is invariant by the action of $W$ (and if $W$ is irreducible, this is unique up to multiplication by positive real numbers). We fix one such form $(\cdot,\cdot)$; if $A$ is a linear transformation of $\hh$ we say $A$ is \emph{Hermitian} if it is equal to its adjoint with respect to this fixed $W$-invariant positive definite Hermitian form (this is actually independent of our choice of form if $W$ is irreducible).

Given a Hermitian operator $A$ on $\hh$, and writing $(x,y)$ for the $W$-invariant Hermitian form fixed above, we obtain a Hermitian form $(Ax,y)$ on $\hh$. If this satisfies
$$(Ax,x) \geq 0 \quad \hbox{for all $x \in \hh$}$$ then we say $A$ is Hermitian \emph{positive semi-definite}. A sum of Hermitian positive semi-definite operators is Hermitian positive semi-definite. 

\subsection{Coxeter numbers of reflection subgroups}

A \emph{reflection subgroup} of $W$ is a subgroup $W' \subseteq W$ that is generated by $R'=W' \cap R$. There is a unique $W'$-stable subspace $\hh' \subseteq \hh$ such that $\hh=\hh^{W'} \oplus \hh'$, and we call $W'$ \emph{irreducible} if $\hh'$ is an irreducible $\C W'$-module. Let $d$ be a positive integer and let $\zeta$ be a primitive $d$th root of unity. We say that $d$ is a \emph{regular number} for $W$ if there is a vector $$v \in \hh^\circ=\hh \setminus \bigcup_{H \in \mathcal{A}} H$$ and an element $w \in W$ with $wv=\zeta v$; this is equivalent to $d$ dividing the same number of degrees of $W$ as codegrees of $W$. In this case any $w \in W$ with $wv=\zeta v$ is called a \emph{regular element of order $d$}, and the regular elements of order $d$ are all conjugate to one another in $W$. A \emph{Coxeter element} of a well-generated complex reflection group is a regular element corresponding to the regular number $d_n$ (where $d_n$ is the unique largest degree of $W$; we note that if W is irreducible then according to the  
classification, $h=d_n$ precisely if $W$ is well-generated, though we know of no reason for this coincidence).

Finally, we recall Springer's result that if $d$ is a positive integer and $\zeta$ is a primitive $d$th root of $1$, then the maximum dimension of the $\zeta$-eigenspace for $w$ ranging over $W$ is the number of degrees of $W$ divisible by $d$. In particular, if $d$ does not divide any degree then there is no $w \in W$ for which $\zeta$ is an eigenvalue.

\begin{theorem}
\begin{enumerate} \item[(a)] Let $W$ be an irreducible complex reflection group with Coxeter number $h$ and let $W' \subseteq W$ be an irreducible reflection subgroup with Coxeter number $h'$. Then $h' \leq h$ with equality if and only if $W'=W$.
\item[(b)] Let $W$ be an irreducible well-generated complex reflection group and let $w \in W$ be a Coxeter element. Then for all reflection subgroups $W' \neq W$ of $W$, $w \notin W'$.
\end{enumerate}
\end{theorem}
\begin{proof}
We first prove (a). Let $\hh' \subseteq \hh$ be the irreducible reflection representation of $W'$. The difference $h-h'$ is the scalar by which
$$\sum_{r \in R \setminus R'} 1-r$$ acts on $\hh'$.  For each reflecting hyperplane $H$ of $W$ let $$W_H=\{r \in W \ | \ r(p)=p \quad \hbox{for all $p \in H$} \}$$ be the (cyclic) group fixing $H$ pointwise, so that every non-trivial element of $W_H$ is a reflection, and let $W_H'=W' \cap W_H$. The element 
$$\sum_{r \in W_H \setminus W_H'} (1-r)$$ is then Hermitian positive semi-definite: upon simultaneously diagonalizing the elements of $W_H$ this reduces to the fact that
$$\sum_{\zeta \in C_n \setminus C_m} (1-\zeta)=n-m$$ where $C_n$ is the group of $n$th roots of unity in $\C^\times$. Hence
$$\sum_{r \in R \setminus R'} (1-r)$$ is Hermitian positive semi-definite. The scalar by $h-h'$ by which it acts on $\hh'$ is therefore non-negative. We suppose this scalar is $0$. Then for all $x \in \hh'$ we have
$$0=\left(\sum_{r \in R \setminus R'} (1-r) x,x \right)=\sum_{H \in \mathcal{A}} \left(\sum_{r \in W_H \setminus W_H'} (1-r)x,x \right)$$ and since for each $H \in \mathcal{H}$ the summand on the right-hand side is non-negative,
$$\left(\sum_{r \in W_H \setminus W_H'} (1-r) x,x \right)=0 \quad \hbox{for all $H \in \mathcal{A}$ and all $x \in \hh'$. }$$ Let $x=x_1v_1+\cdots+x_n v_n$ where $v_2,\dots,v_n$ are a basis of the fix space $H$ of $W_H$ and $v_1$ is a common eigenvector for $W_H$, normalized so that $(v_i,v_j)=\delta_{ij}$ is the Kronecker delta. The previous equation implies $x_1=0$ if $W_H' \neq W_H$, or in other words $x \in H$ for all such $H$. Hence $\hh' \subseteq H$ if $W_H' \neq W_H$. Thus $\hh'$ is fixed by $W_H$ for all $H$ with $W_H' \neq W_H$. But by definition $\hh'$ is $W_H$-stable if $W_H'=W_H$, so that $\hh'$ is $W$-stable, which by irreducibility of $\hh$ gives $\hh'=\hh$. Thus $\hh'$ is not contained in \emph{any} $H \in \mathcal{A}$, and from the previous reasoning we must have $W_H=W_H'$ for all $H$, or in other words $W'=W$.

For (b), we observe that for an arbitrary (not necessarily irreducible) reflection subgroup $W'$ of $W$, the degrees of $W'$ are the degrees of its irreducible factors, and all of these are at most the respective Coxeter numbers (by examining the classification). These Coxeter numbers are, for proper $W'$, strictly less that the Coxeter number of $W$. The result follows. \end{proof}

\subsection{The groups $G(\ell,m,n)$} Fix $\zeta$ a primitive $\ell$th root of $1$, and for $\mu \in (\Z / \ell \Z)^n$ write $\zeta^\mu$ for the diagonal matrix whose diagonal entries are the numbers $\zeta^{\mu_i}$ for $i=1,2,\dots,n$. Given positive integers $\ell,m$ and $n$ with $m$ dividing $\ell$, the group $G(\ell,m,n)$ consists of all products $\zeta^\mu w$ where $w \in S_n$ is a permutation matrix and $\mu \in (\Z / \ell \Z)^n$ satisfies
$$(\ell/m) \cdot (\mu_1+\dots+\mu_n)=0 \quad \mathrm{mod} \ \ell.$$ Putting $\mu=(m-1,1,0,\dots,0)$ one verifies that the element
$$w=\zeta^\mu (12 \cdots n) \in G(\ell,m,n)$$ does not belong to any reflection subgroup of $G(\ell,m,n)$.

\subsection{The group $G_{15}$}

We refer to chapter 6 of Lehrer and Taylor's book \cite{LT} for the unproved assertions about rank two groups which follow.

The complex reflection group $G_{15}$ has three maximal reflection subgroups (all of which are normal): $G_7$, $G_{13}$, and $G_{14}$, of orders $144$, $96$, and $144$, respectively. The intersection of $G_{14}$ and $G_{13}$ is $G_{12}$, of order $48$, and the intersection of $G_{14}$ and $G_7$ is $G_5$, of order $72$. It follows that the union $G_{13} \cup G_{14}$ has $96+144-48=192$ elements, while the union $G_7 \cup (G_{13} \cup G_{14})$ has at most $144+192-72=264$ elements. Thus there is an element of $G_{15}$ not contained in any reflection subgroup.

We note that the degrees and codegrees for $G_{15}$ are $(12,24)$ and $(0,24)$. So the regular numbers are the divisors of $12$. On the other hand, every divisor of $12$ divides \emph{both} degrees, and hence a regular element is just a scalar matrix. Thus the regular elements in $G_{15}$ are precisely the scalar matrices whose entries are $12$th roots of $1$. These all belong to the subgroup $G_7$ of $G_{15}$, and therefore there is no regular element of $G_{15}$ that does not belong to any proper reflection subgroup.

\subsection{Completion of the proof} Examining the list of complex reflection groups and the data in \cite{Tay} shows that for each irreducible exceptional group $W$, either (i) the group is well-generated, or (ii) there is an integer $d$ that divides more degrees of $W$ than of any proper reflection subgroup, or (iii) $W$ is the group $G_{15}$. For instance, the non-well-generated group $G_{31}$ (also sometimes known as $W(O_4)$) has degrees $8,12,20$, and $24$. According to Taylor \cite{Tay}, each maximal reflection subgroup of $W(O_4)$ is one of $G(4,2,4)$, $W(F_4)$, or $W(N_4)=G_{29}$. But $24$ does not divide any of the degrees of these groups.

In the first two cases our result and Springer's theory of eigenvalues shows that there is an element of $W$ that does not belong to any proper reflection subgroup. By the classification theorem and reduction to the case $W$ irreducible, this proves:
\begin{theorem} \label{appendixthm}
Let $W$ be a finite complex reflection group. Then there is some $w \in W$ such that for each reflection subgroup $W' \neq W$, $w \notin W'$.
\end{theorem}

Note that the second condition required by Proposition \ref{prop fixed point} is immediate. Indeed, an element with a trivial eigenvalue is by definition contained in a proper parabolic subgroup and Steinberg's Theorem states that these are reflection subgroups.

\end{appendix}

\end{document}